\documentclass{jas}
\usepackage[all]{hypcap}
\usepackage{csquotes,tikz,lipsum}
\usetikzlibrary{shapes}

\title{The applicability of equal area partitions of the unit sphere}
\Shorttitle{Applicability of equal area sphere partitions}

\author[1,2,$\star$]{Paul Leopardi,\orcidlink{0000-0003-2891-5969}}
\Shortauthors{P. Leopardi}

\affil[1]{ACCESS-NRI}
\affil[2]{The Australian National University}

\Email{paul.leopardi@anu.edu.au}

\Subdate{23/04/2024}

\Volume{W}
\Issue{X}
\Year{YYYY}
\Artnum{Z}
\DOI{0000.000000}
\Pubdate{DD/MM/YYYY}
\Comby{}

\newcommand{\Real}{\mathbb R}
\newcommand{\Complex}{\mathbb C}
\newcommand{\Sphere}{\mathbb S}

\newcommand{\Dim}{d}

\newcommand{\Energy}{\operatorname{E}}
\newcommand{\EQCode}{\operatorname{EQP}}
\newcommand{\norm}[1]{\left\| #1\right\|}

\newcommand{\Partition}{\operatorname{EQ}}

\begin{document}

	\maketitle

	\begin{abstract}
This paper addresses the idea of the applicability of mathematics,
using, as a case study, a construction and software package that partition the unit sphere into regions of equal area.
The paper assesses the applicability of this construction and software by examining citing works,
including papers, dissertations and software.
	\end{abstract}


\section{Introduction}
\paragraph*{The unreasonable effectiveness of mathematics?}
A well known paper by Wigner \cite{wigner1960unreasonable}
claims that there is something mysterious or even miraculous about the appropriateness of mathematics as a
language for describing nature \ldots
\begin{quote}
The miracle of the appropriateness of the language of mathematics for the formulation of
the laws of physics is a wonderful gift which we neither understand nor deserve.
We should be grateful for it and hope that it will remain valid in future research and
that it will extend, for better or worse,
to our pleasure even though perhaps also to our bafflement,
to wide branches of learning.
\cite{wigner1960unreasonable}
\end{quote}

Hamming \cite{hamming1980unreasonable} make some attempts at explanation,
but comes to a similar conclusion to Wigner \ldots
\begin{quote}
Some partial explanations \dots
\begin{enumerate}
 \item We see what we look for.
 \item We select the kind of mathematics to use.
 \item Science in fact answers comparatively few problems.
 \item The evolution of man provided the model.
\end{enumerate}
\ldots
From all this I am forced to conclude both that mathematics is unreasonably effective and that
all of the explanations I have given when added together simply are not enough to explain
what I set out to account for. \ldots
The logical side of the nature of the universe requires further exploration.
\cite{hamming1980unreasonable}
\end{quote}

On the other hand, Arnold \cite{arnold2014teaching} and Borovik \cite{borovik2021mathematician}
quote Gelfand's assertion that mathematics is unreasonably \emph{in}effective in biology.

Arnold \cite{arnold2014teaching}:
\begin{quote}
Here we can add a remark by I.M. Gel'fand:
there exists yet another phenomenon which is comparable in its inconceivability with
the inconceivable effectiveness of mathematics in physics noted by Wigner --
this is the equally inconceivable ineffectiveness of mathematics in biology.
\end{quote}

Borovik: \cite{borovik2021mathematician}:
\begin{quote}
This paper is an attempt to answer the question
\begin{quote}
Should we accept Israel Gelfand's assessment of the role of mathematics in biology?
\end{quote}
And my answer is
\begin{quote}
Yes, we should, for the time being: mathematics is still too weak for playing in biology the role it ought to play.
\end{quote}
\end{quote}

So what makes mathematics effective or ineffective in applications, and is the effectiveness or ineffectiveness reasonable or unreasonable?
\paragraph*{A more pragmatic approach.}
The practice of applied mathematics usually takes a more pragmatic approach, especially when dealing with models and approximations.
In scientific modelling, it is often stated that ``All models are wrong but some are useful'' \cite{box1979robustness},
meaning that a completely faithful model may be unattainable, but it may be possible to build a parsimonious model that reflects
the key phenomena or most important aspects of the system being modelled.

Much can and has been said about the construction of models of systems that are both
fit for purpose and mathematically tractable \cite{lisciandra2017robustness}. Modelling often involves approximation, in the
sense of neglecting some aspects of the systems, and idealization, that is making strictly incorrect
assumptions that still preserve the important aspects to be understood \cite{macleod2021applicability,norton2012approximation,Row2011approximations}.
The processes of approximation and idealization may then result in a model that can be described by
known or at least constructable mathematics. The final step for a predictive model would then be to ensure that the mathematical formulation is tractable, in the sense that it results in a reasonable trade-off between computational effort and accuracy \cite{aaronson2013philosophers, azzouni2000applying,rice1976algorithm,traub1991information}.

Models also often involve approximations in the mathematical sense of the word.
The approximation of functions from noisy and incomplete data \cite{aste2015techniques},
and the approximate solution of underdetermined systems of equations \cite{jokar2008exact},
including the solution of differential or integral equations with noisy and incomplete initial data, has long been a subject of study in
statistics, applied mathematics, and machine learning. The subject of approximation theory deals with the best approximation within a function space \cite{trefethen2019approximation}, and the theory of information-based complexity explicitly deals with the inherent trade-off between the cost of function evaluation versus the error in approximation \cite{traub1991information,wozniakowski2009information}.
It is therefore usually the case that a new approximation method, algorithm or construction is accompanied by an analysis of its applicability to known abstract problems, often in terms of its asymptotic cost versus rate of convergence with respect to approximation error within a known setting (e.g. \cite{brauchart2015covering}).

\paragraph*{A case study.}
This paper examines the applicability of a geometric construction: an equal area partition of a unit higher-dimensional sphere, and its associated distribution of points. The construction is described in a 2006 paper published in Electronic Transactions on Numerical Analysis \cite{leopardi2006partition}, and is analyzed in more detail in a PhD thesis of 2007 \cite{leopardi2007distributing}, a paper of 2009
\cite{leopardi2009diameter} and follow-up papers \cite{leopardi2013discrepancy,leopardi2017optimal}.
As at 20 August 2024, the 2006 paper \cite{leopardi2006partition} has 346 citations listed in Google Scholar,
224 citations in Scopus, 179 in Web of Science, and 43 in MathSciNet.

Citations appear in
Geophysical Journal International,
Global Change Biology,
IEEE Transactions on Audio Speech and Language Processing,
Journal of Approximation Theory,
Journal of the Atmospheric Sciences,
Journal of Computational Chemistry,
Journal of Differential Equations,
Mathematics of Computation,
Radio Science,
RNA Journal, and elsewhere.

Note: some of the papers described below are accompanied by abbreviated 20 August 2024 citation counts of the form (G: g, S: s, W: w, M: m), for the Google Scholar, Scopus, Web of Science, and MathSciNet counts respectively. For example, the abbreviated counts for the 2006 paper \cite{leopardi2006partition} are (G: 346, S: 224, W: 179, M: 43).

The citations of the 2006 paper, the 2007 thesis, and the follow-up papers are generally of three types:
\begin{enumerate}
 \item Application of the constructions to specific problems;
 \item Evaluation of the constructions described by the paper, including comparisons with other constructions; and
 \item Passing mention of the paper, sometimes with a short description.
\end{enumerate}
This paper is mostly concerned with the first two types of citations.

\section{Preliminaries}
\subsection{Some related problems}\label{sec:related-problems}
The problem of finding an equally distributed set of $N$ points on a circle is solved easily: just use points arranged uniformly at
an angle of $\frac{N}{2\pi}$.
In contrast, on a unit sphere $\Sphere^d \subset \Real^{d+1}$ with $d>1$, not only is the problem harder to solve,
it is harder to pose precisely.
There are a number of related problems, each of which gives rise to a different sense of equal distribution \cite{borodachov2019discrete,SafK97}.
These problems are often stated in terms of sequences of \emph{spherical codes}, where each spherical code is a finite set of $N$ points on
the unit sphere $\Sphere^d$, and we are often interested in some asymptotic value related to each sequence as $N$ goes to infinity.

\begin{figure}[!ht]
\begin{center}
\includegraphics[angle=0,width=80mm]{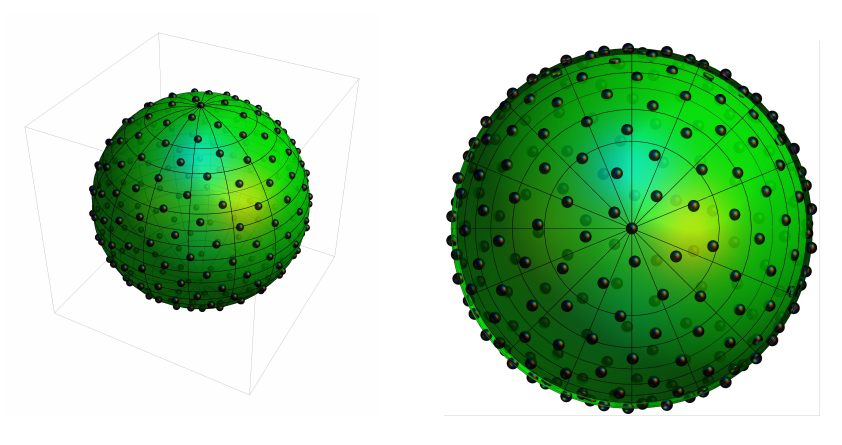}
\caption{Brauchart, ``Spherical Fibonacci lattices,''
\cite{brauchart2012low}.}
\end{center}
\label{fig:brauchart2012low}
\end{figure}

\begin{figure}[!ht]
\begin{center}
\includegraphics[angle=0,width=120mm]{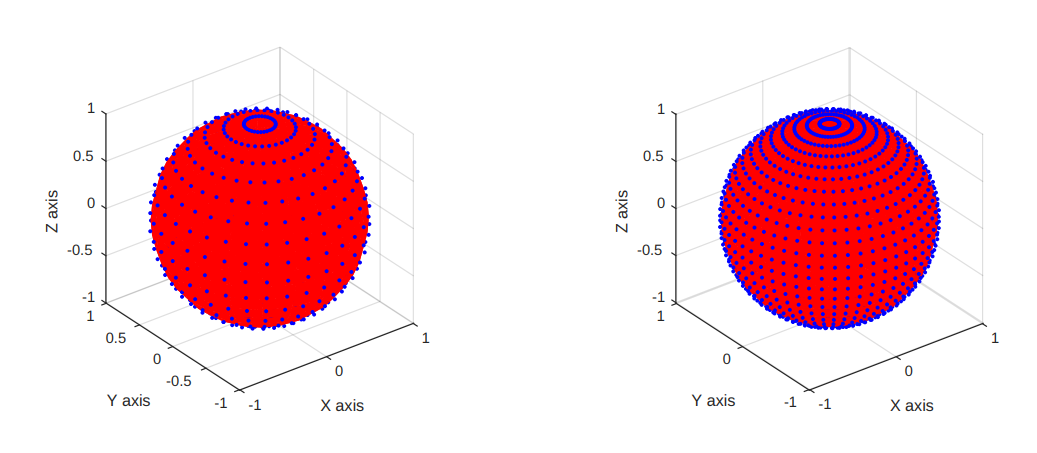}
\caption{Themistoclakis and Van Barel, ``Examples of the tensor product Gauss–Legendre quadrature nodes related to
degrees of precision 31 and 51, i.e., having $N = 512$ (left) and $N = 1352$ (right) points,''
\cite{themistoclakis2018optimal}.}
\end{center}
\label{fig:themistoclakis2018optimal}
\end{figure}

\begin{figure}[!ht]
\begin{center}
\includegraphics[angle=0,width=80mm]{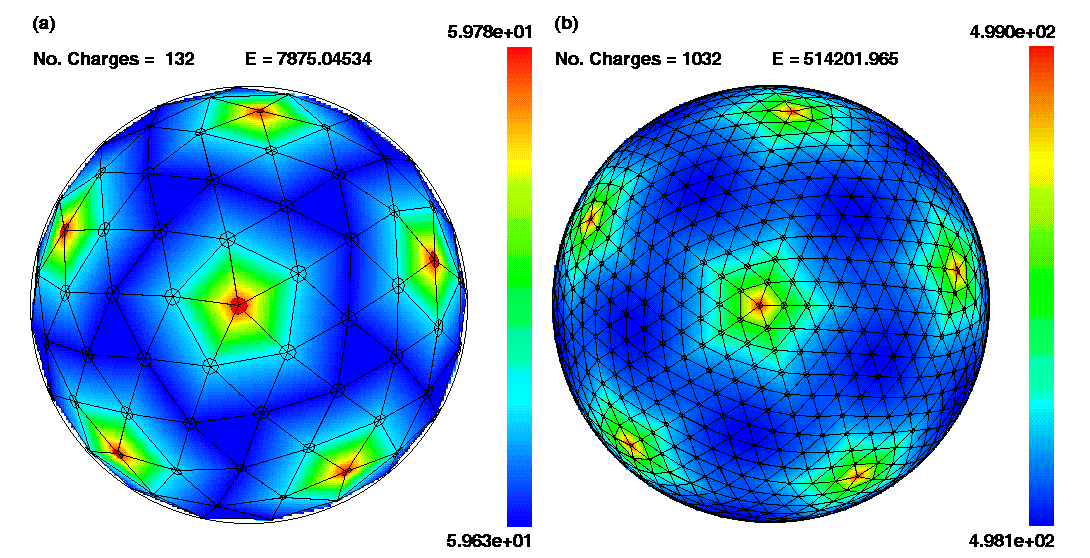}
\caption{Atschuler et al. ``Lattice configurations for 132 (a) and 1032 (b) charges,''
\cite{altschuler1997possible}.}
\end{center}
\label{fig:altschuler1997possible}
\end{figure}

\begin{figure}[!ht]
\begin{center}
\includegraphics[angle=0,width=50mm]{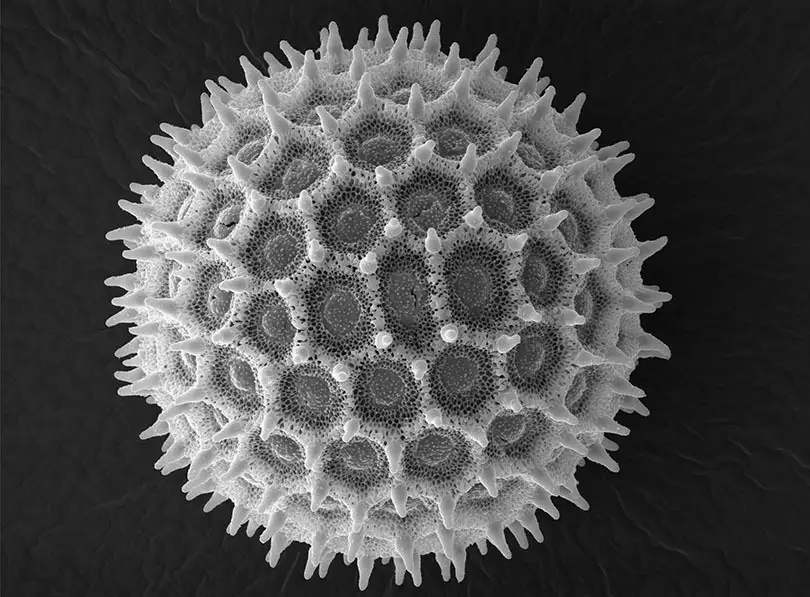}
\caption{Dartmouth College Electron Microscope Facility, ``A grain of pollen from Morning Glory flowers,''
\cite{dartmouth2024electron}.}
\end{center}
\label{fig:dartmout2024electron}
\end{figure}

\paragraph*{Asymptotic equidistribution in measure and related quadrature problems.}
Find a sequence of spherical codes that converges to the uniform distribution in the sense of convergence
in measure
(Figure \ref{fig:brauchart2012low}).
Such convergence is usually defined in terms of a discrepancy such as the spherical cap discrepancy \cite{Bra01,GraT93,wong1997sampling}.
If the spherical cap discrepancy of the spherical codes in the sequence converges to zero, then the sequence is said to be
\emph{asymptotically uniformly distributed} \cite{SafK97}, or \emph{asymptotically equidistributed} \cite{DamG03,MarO09},
or \emph{weak-star convergent} \cite[Definition 2.11.3]{leopardi2007distributing}.
\paragraph*{Interpolation and related function approximation problems.}
Find a sequence of spherical codes such that a function approximation within a certain function space converges at a
certain rate.
Examples are Lagrange polynomial interpolation within continuous function spaces \cite{AnCSW10,bos1999limiting,WomS01}
and least squares approximation within $L_{\infty}$ \cite{JetSW98b,JetSW99,sloan2000constructive,themistoclakis2018optimal} (Figure \ref{fig:themistoclakis2018optimal}).
\paragraph*{The Thomson and related energy minimization problems.}
Minimize the energy of $N$ equally charged particles on a sphere, with respect to some potential
(Figure \ref{fig:altschuler1997possible}).
\paragraph*{The Tammes problem and packing of spherical caps.}
Given a fixed radius, how many non-overlapping spherical caps with that radius can be placed onto a unit sphere
(Figure \ref{fig:dartmout2024electron})?
This radius is called the \emph{packing radius} of the spherical code formed by the centres of the caps.

\begin{figure}[!ht]
\begin{center}
\includegraphics[angle=0,width=40mm]{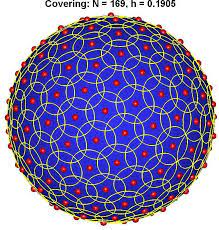}
\caption{Saff and Womersley, ``Covering of a sphere with 169 equal spherical caps,''
\cite{saff2013optimal}.}
\end{center}
\label{fig:saff2013optimal}
\end{figure}

\paragraph*{The covering problem.}
Given a fixed radius, how few overlapping spherical caps with that radius are needed to cover a unit sphere
(Figure \ref{fig:saff2013optimal})?
This radius is called the \emph{covering radius} or \emph{mesh norm} of the spherical code formed by the centres of the caps.

\subsection{Some history}
The history of constructions aimed at solving the problems posed in Section \ref{sec:related-problems} is quite involved.
See also the 2019 book by Borodachov et al. \cite[Chapters 6 and 7]{borodachov2019discrete}.
\paragraph*{Equidistribution without separation.}

\begin{figure}[!ht]
\begin{center}
\includegraphics[angle=0,width=60mm]{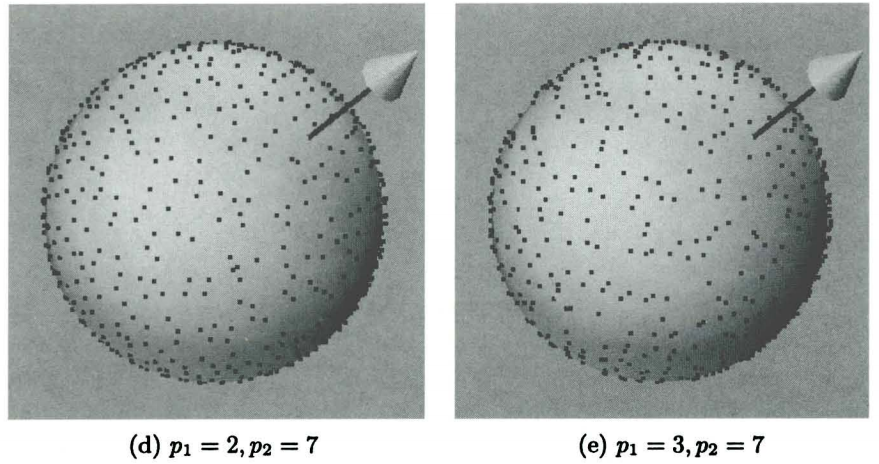}
\caption{Wong and colleagues ``Halton points with different bases on the sphere ($n = 1000$),'' \cite{wong1997sampling}.}
\end{center}
\label{fig:wong1997sampling}
\end{figure}

Many constructions for $\Sphere^2$ yield an asymptotic equidistribution,
e.g. Hammersley, Halton, $(t,s)$ etc. sequences mapped to the sphere  \cite{wong1997sampling} (Figure \ref{fig:wong1997sampling}).
\paragraph*{Separation without equidistribution.}

\begin{figure}[!ht]
\begin{center}
\includegraphics[angle=0,width=35mm]{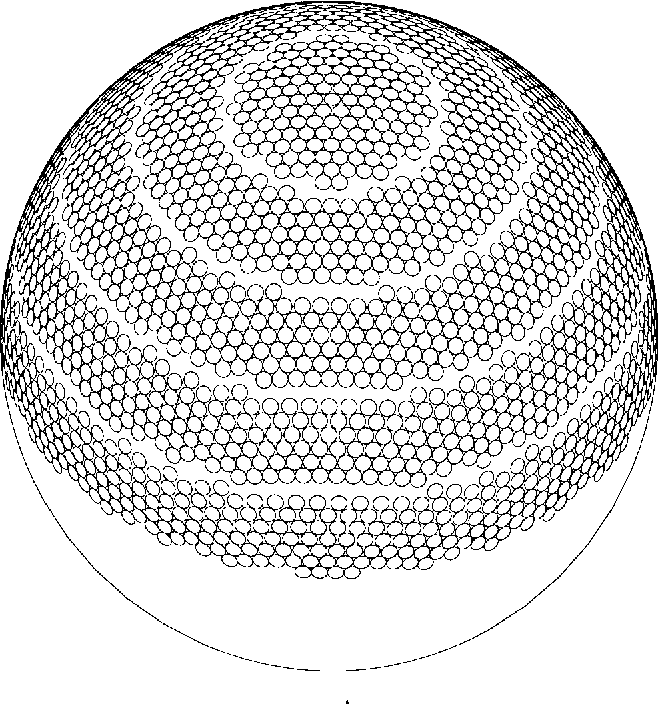}
\caption{Hamkins and Zeger, ``The wrapped spherical code $\mathcal{C}_{W}^{\Lambda_2}$ $(3; 0:05)$,'' \cite{HamZ97a}.}
\end{center}
\label{fig:HamZ97a}
\end{figure}

Hamkins \cite{Ham96} and Hamkins and Zeger \cite{HamZ97a} constructed $\Sphere^{\Dim}$ codes with asymptotically optimal packing density
(Figure \ref{fig:HamZ97a}).

\paragraph*{Equal area partitions.}

Alexander \cite{Ale72} asserts the existence of a diameter bounded set of equal area partitions of $\Sphere^2$,
sketching a construction based on the cubed sphere.

Stolarsky \cite{Sto73}, Beck and Chen \cite{BecC87} and Bourgain and Lindenstrauss \cite{BouL88}
each go on to assert the existence of a diameter bounded set of equal area partitions of $\Sphere^{\Dim}$ without giving an explicit construction.

\begin{figure}[!ht]
\begin{center}
\includegraphics[angle=0,width=50mm]{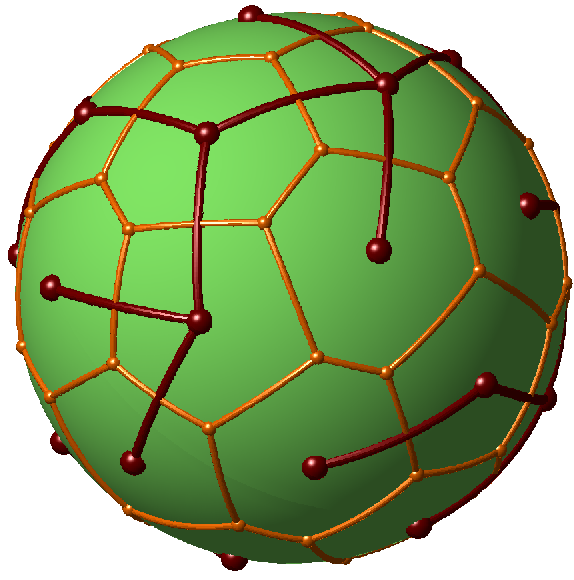}
\caption{Leopardi. ``Step 5 of the Feige-Schechtman construction'' \cite{FeiS02,leopardi2007distributing}.}
\end{center}
\label{fig:FeiS02}
\end{figure}

Feige and Schechtman \cite{FeiS02} describe a construction as part of
an argument about the optimality of a solution of the Max-Cut problem in graph theory that can be modified into a
construction of a diameter bounded set of equal area partitions of $\Sphere^{\Dim}$ (Figure \ref{fig:FeiS02}).

\begin{figure}[!ht]
\begin{center}
\includegraphics[angle=0,width=60mm]{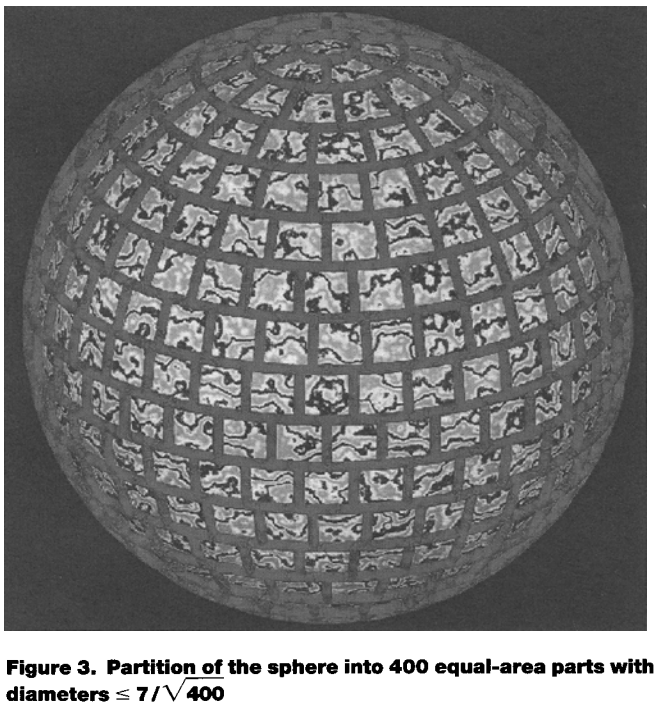}
\caption{Saff and Kuijlaars. ``Partition of the sphere into 400 equal-area parts with diameters $\leq 7 / \sqrt{400}$''
\cite{SafK97}.}
\end{center}
\label{fig:SafK97}
\end{figure}

\begin{figure}[!ht]
\begin{center}
\includegraphics[angle=0,width=60mm]{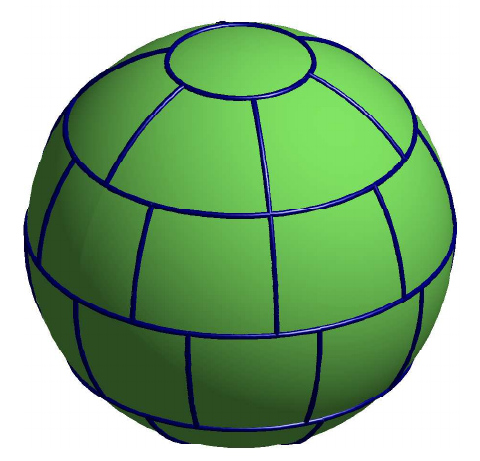}
\caption{``Partition $\Partition(2,33)$''
\cite{leopardi2006partition}.}
\end{center}
\label{fig:leopardi2006partition}
\end{figure}

The $\Partition(d, N)$ recursive zonal
partition of the sphere $\Sphere^d$ into $N$ regions of equal area described in \cite{leopardi2006partition} and analyzed in \cite{leopardi2007distributing,leopardi2009diameter} is based on Zhou's 1994 construction for $\Sphere^2$ \cite{RakSZ94,Zho95} (Figure \ref{fig:SafK97})
as modified by Saff, and Sloan's sketch of a partition of $\Sphere^3$ \cite{Slo03} (Figure \ref{fig:leopardi2006partition}).
%

The 2006 paper \cite{leopardi2006partition} describes the construction of the $\Partition(d, N)$ partition.
The paper also provides estimates and numerical examples of the maximum diameter of each region in each partition of $\Sphere^d$ into $N$ regions for $N \leq 100\ 000$ for $d=2,3,4$, and also for $N=2^k$ for $k=1 \dots 10$ and $d=1 \ldots 8$.
The maximum diameter is a good estimate for twice the covering radius.
\paragraph{The partition algorithm.}

The recursive zonal equal area partition algorithm is recursive in dimension $d$.
For $d > 1$ it uses the idea of a ``collar'' -- an annulus on the sphere arranged symmetrically about the North-South polar axis.

The 2006 paper \cite{leopardi2006partition} provides a detailed description of the partition algorithm $\Partition(d,N)$, but a brief pseudocode description is
\small{}
\begin{quote}
\noindent
$\mathbf{if}$ $N = 1$ $\mathbf{then}$
\begin{verse}
There is a single region which is the whole sphere;
\end{verse}
$\mathbf{else}$ $\mathbf{if}$ $d = 1$ $\mathbf{then}$
\begin{quote}
Divide the circle into $N$ equal segments;
\end{quote}
$\mathbf{else}$
\begin{quote}
Divide the sphere into zones,
each the same area as an integer number of regions:
\begin{enumerate}
 \item
Determine the colatitudes of polar caps,
 \item
Determine an ideal collar angle,
 \item
Determine an ideal number of collars,
 \item
Determine the actual number of collars,
 \item
Create a list of the ideal number of regions in each collar,
 \item
Create a list of the actual number of regions in each collar,
 \item
Create a list of colatitudes of each zone;
\end{enumerate}
Partition each spherical collar into regions of equal area,
using the recursive zonal equal area partition algorithm for dimension $d-1$;
\end{quote}
$\mathbf{end if}$.
\end{quote}
\normalsize{}

\begin{figure}[!ht]
\begin{center}
\includegraphics[width=120mm]{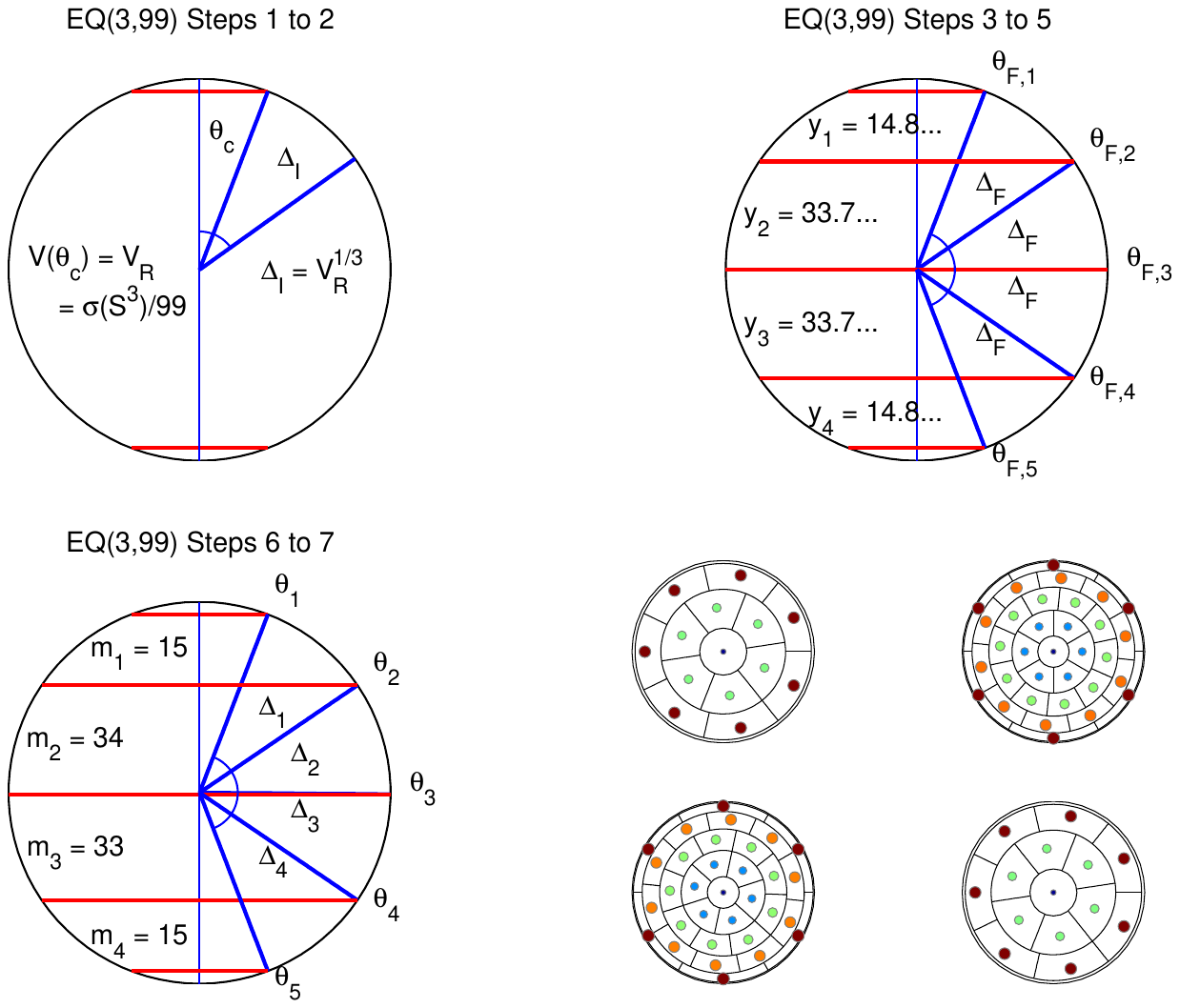}
\end{center}
\caption{Partition algorithm for $\Partition(3,99)$}
\label{algorithm-3-99}
\end{figure}
\begin{figure}[!ht]
\begin{center}
\includegraphics[angle=270,width=100mm]{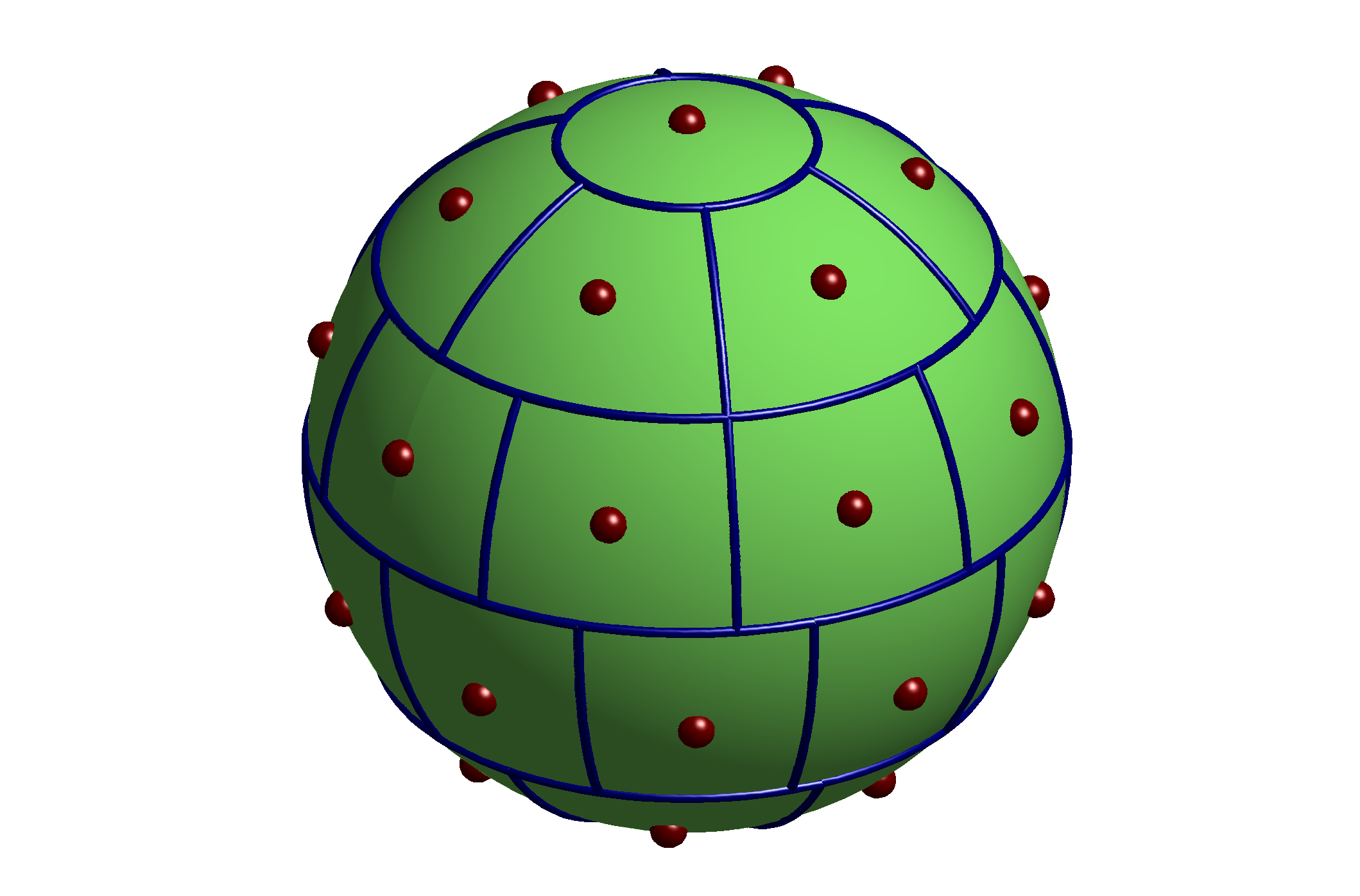}
\caption{``EQ code $\EQCode(2, 33)$, showing the partition $\Partition(2, 33)$''
\cite{leopardi2007distributing}.}
\end{center}
\label{fig:leopardi2007distributing}
\end{figure}

Figure \ref{algorithm-3-99} is an illustration of the algorithm for $\Partition(3,99)$,
with step numbers corresponding to the step numbers in the pseudocode.
\paragraph{Spherical codes from equal area partitions.}
The 2007 thesis \cite{leopardi2007distributing} describes the partition in more detail, describes the spherical codes
$\EQCode(d, N)$ consisting of a central point of each region of $\Partition(d, N)$, and proves that the sequences of these
codes are asymptotically equidistributed for each $d$ \cite[Theorem 5.4.1]{leopardi2007distributing} (Figure \ref{fig:leopardi2007distributing}). Despite being unpublished, as at 20 August 2024 the thesis has 78 citations on Google Scholar.

Chapters 3 to 5 of the 2007 thesis also include the following statements, estimates and numerical examples:

Chapter 3 contains estimates of the maximum diameter of each region, as per \cite{leopardi2006partition}, with proofs.
Section 3.10 lists numerical results on the maximum diameter of regions, as per \cite{leopardi2006partition}.

Section 4.2 includes an explanation of why the spherical codes are \emph{not good} for Lagrange polynomial interpolation, in terms of the condition number of the Gram matrix corresponding to each polynomial degree (see also \cite{WomS01}).
 The section remarks that the Gram matrix is often singular to machine precision,
 and also includes a statement and proof of the minimum polynomial degree for which Lagrange interpolation must fail.
Section 4.3 estimates the minimum distance between code points, which is also a good estimate of twice the packing radius.
 Numerical examples are given for $d=2,3,4$ and $N$ up to $20\ 000$.
Section 4.4 gives estimates of the packing density, which is also related to the packing radius. This section also includes numerical examples for $d=2,3,4$ and $N$ up to $20\ 000$.

Section 5.4 includes estimates of the spherical cap discrepancy and the Riesz energy, with numerical examples of the Riesz $d-1$ energy for $d=2,3,4$ and $N$ up to $20\ 000$.
Here the Riesz $s$-energy of a finite set $\mathcal{X} \subset \Sphere^d$ of size $N$ is defined as
\begin{align*}
\Energy_s(\mathcal{X})
&=
\frac{1}{N^2}
\sum_{x \in \mathcal{X}} \mathop{\sum_{y \in \mathcal{X}}}_{y\neq x}
\norm{x - y}^{-s},
\end{align*}
using the usual Euclidean norm on $\Real^{d+1}$.

\paragraph*{Matlab code.}

The Recursive Zonal Equal Area Sphere Partitioning Toolbox \cite{leopardi2005recursive,leopardi2017recursive,leopardi2024recursive}
is a Matlab toolbox that was released in 2005 to accompany the paper \cite{leopardi2006partition} and PhD thesis \cite{leopardi2007distributing}. The earlier history of the code, including the original Maple prototype, can be seen in the \texttt{CHANGELOG} file \cite{leopardi2024recursive}. The following remarks refer to the situation at 20 August 2024.

Google Scholar lists four citations to the toolbox, excluding self-citations \cite{deboy2010acoustic,deboy2010acousticdiploma,rajan2023analysis,vianello2018global}.
The SourceForge URL of the toolbox \cite{leopardi2005recursive} is mentioned in 17 other theses and papers indexed by Google Scholar \cite{etayo2021spherical,hardin2016comparison,jenkins2012local,kurz2017linear,lin2023distributed,lin2024sketching,mozdzynski2007new,pomberger2014estimating,prats2017configuracions,rajan2023analysis,schouten2021simulation,smale2007diffusion,szwajcowski2021error,van2012directional,wu2018four,wu2019four,wu2020hybrid}, excluding self-citations.
Of these theses and papers only six
\cite{kurz2017linear,pomberger2014estimating,prats2017configuracions,rajan2023analysis,smale2007diffusion,szwajcowski2021error} contain an attributed citation to the toolbox in their References section.
Interestingly, one thesis \cite{wu2019four} and two papers \cite{wu2018four,wu2020hybrid} mistakenly call the toolbox ``EASP'' and do not cite the author.

The GitHub project for the toolbox \cite{leopardi2024recursive} has 8 forks in GitHub, and is mentioned in one other paper indexed by Google Scholar \cite{schouten2021simulation}.
Code from the toolbox is also included in at least 10 other GitHub projects \cite{Frey2016,Ivanov2021,Kam2018,Kunc2020,Kurz2023,Li2021,Na2024,Pan2020,Viard2020,Wang2023}.
 Unfortunately, SourceForge does not support code searches across its repositories, so the number of SourceForge projects that include code from the toolbox remains unknown. The same is true for GitLab public repository hosting.

The toolbox is also mentioned in the documentation for FERUM \cite{bourinet2010ferum} but without attribution.
 The Matlab source code for FERUM 4.1 \cite{Bourinet2009ICOSSAR} contains a subset of the toolbox code.
 There is also a copy of the toolbox code at the Lamont-Doherty Earth Observatory \texttt{clifford.ldeo.columbia.edu} web site \cite{Steiger2021}.

\subsection{Follow-up papers and generalizations}

The 2009 paper \cite{leopardi2009diameter}, based on the 2007 thesis \cite{leopardi2007distributing},
proves diameter bounds for both the $\Partition(d, N)$ sphere partition described in the 2006 paper \cite{leopardi2006partition},
and a modified version of the construction of Feige and Schechtman as described in the thesis. Citations: (G: 34, S: 0, W: 13, M: 11).

A 2013 paper \cite{leopardi2013discrepancy}, following the arguments in Chapter 5 of the 2007 thesis \cite{leopardi2007distributing},
shows that a sequence of spherical codes with a well behaved upper bound on discrepancy and a well behaved lower bound on separation,
such as the sequence of $\EQCode(d, N)$  codes, satisfies an upper bound on the Riesz $s$-energy. Citations: (G: 21, S: 13, W: 12, M: 10).

A second 2013 paper \cite{leopardi2013dolomites}
generalizes the paper \cite{leopardi2013discrepancy} in the sense that it proves that,
for a smooth compact connected $d$-dimensional Riemannian manifold $M$, if $0 \leq s \leq d$ then an asymptotically equidistributed sequence of finite subsets of $M$ that is also well-separated yields a sequence of Riesz $s$-energies that converges to the energy double integral. In this case, the Riesz $s$-energy is
defined using the geodesic distance on $M$.
Citations: (G: 3, S: 0, W: 0, M: 0).

A 2017 joint paper with Gigante \cite{gigante2017diameter}
generalizes the partition results of  \cite{leopardi2007distributing,leopardi2009diameter}
by combining the Feige and Schechtman construction with David's and Christ's dyadic cubes to yield a partition algorithm for
connected Ahlfors regular metric measure spaces of finite measure. Citations: (G: 33, S: 19, W: 18, M: 18).

A second 2017 joint paper with Sommariva and Vianello \cite{leopardi2017optimal}
proves that good covering point configurations on the 2-sphere are optimal polynomial meshes,
and extracts Caratheodory-Tchakaloff submeshes for compressed least squares fitting.
This implies that the point sets generated by the construction of the 2006 paper \cite{leopardi2006partition}
are optimal polynomial meshes. The paper also provides numerical examples where submeshes based on these point sets are used to construct positive weight quadrature rules. Citations: (G: 6, S: 4, W: 5, M: 1).
\section{Evaluations and improvements}
\subsection{Evaluations}\label{sec:evaluations}
Many of papers citing the 2006 paper \cite{leopardi2006partition} and its related papers conduct one of two types of evaluation:
\begin{enumerate}
 \item They evaluate methods that use the constructions described in \cite{leopardi2006partition} against one or more completely different methods as they apply to the problem being solved in the paper. This type of evaluation is most frequently seen in applications oriented papers and is treated in the Section \ref{sec:applications}.
 \item They evaluate the constructions described in \cite{leopardi2006partition} against similar constructions, especially
 in relation to one or more of the related problems listed in Section \ref{sec:related-problems}.
 Some examples of this type of evaluation follow.
\end{enumerate}

A 2009 paper by Marantis and colleagues \cite{marantis2009comparison} compares three different
point distributions on $\Sphere^2$, including $\EQCode(2,N)$, by using test samples of 240 points and using them to reconstruct a function defined by spherical harmonics up to degree 8: ``it is subsequently sampled with the three proposed sample point distributions and the pattern is reconstructed using the estimated harmonic coefficients.''
Unfortunately the paper does not explicitly state the reconstruction method used.
The $\EQCode(2,240)$ reconstruction fails badly. The paper makes an attempt to explain this.

A 2016 paper by Rachinger and colleagues \cite{rachinger2016phase} compares different
``constellations'' of points on hyperspheres in a complex vector space.
The case described in the paper is a hypersphere in $\Complex^3$,
equivalent to the sphere $\Sphere^5 \subset \Real^6$.
The $\EQCode(5,64)$ and $\EQCode(5,512)$ codes are compared to constellations
obtained via $k$-means clustering, potential minimization, and per-antenna phase shift keying. Curiously, the paper calls the the $\EQCode$ codes ``EQPA constellations.''
The constellations are compared in terms of construction complexity,
capacity, minimum distance, and power efficiency: ``\ldots EQPA works in such a way that the distribution of points becomes more and more uniform as the constellation size increases. This algorithm profits from packing the hypersphere more densely.''

A related 2016 paper by Sedaghat and colleagues \cite{sedaghat2016continuous} compares the $\EQCode$ codes to codes created by spherical K-means clustering with respect to performance of a wireless communication scheme called Phase Modulation on the hypersphere: ``\ldots  the codes obtained by the spherical K-means algorithm have much better
performance than the EQ codes. Note that EQ codes have the advantage that they can be constructed much more easily than
K-means codes.''

One of the most comprehensive comparisons of constructions for spherical codes on $\Sphere^2$ is found in the 2016 paper of Hardin, Michaels and Saff \cite{hardin2016comparison},
which examines quadrature, energy, packing and covering properties of a number of such constructions.
The paper shows that the EQ point sets generated by the construction described in \cite{leopardi2006partition} are not
only equidistributed and well-separated, but they also perform well with respect to energy, with numerical behaviour comparable to
empirically optimal point sets. For the logarithmic and Coulomb potentials, ``the generalized spiral and zonal equal area points perform the best of the algorithmically generated points.'' For the Riesz $s$-energy with $s=2$, ``the generalized spiral, zonal equal area, and equal area icosahedral points perform the best.'' For the Riesz $s$-energy with $s=3$,
``the equal area icosahedral points outperform the spiral and zonal equal area points of the algorithmically generated
configurations. This is expected because their Voronoi decomposition is closest to the regular hexagonal lattice.''
For more detailed proofs, see the 2017 PhD thesis of Michaels \cite{michaels2017node}.
For more context, see the 2019 book by Borodachov and colleagues \cite[Chapter 7]{borodachov2019discrete}.

\subsection{Improvements}

For the sphere $\Sphere^2$ the \emph{diamond ensemble} \cite{beltran2020diamond} is a construction for spherical codes resembling the $\EQCode(2,N)$ codes, where the code is constructed
directly and not via an equal area partition.
Similarly to Zhou's construction \cite{RakSZ94,Zho95} and the $\EQCode(2,N)$ codes, the code points are distributed amongst the north and south poles and a small number of parallels of latitude. On each parallel, the code points are equally spaced.
Unlike Zhou's construction and the $\EQCode(2,N)$ codes, the code points on each parallel are offset by a random angle, and the number of code points per parallel are chosen to minimize the expected logarithmic energy. The diamond ensemble can then
be used to construct an equal area partition similar to Zhou's construction or
the $\Partition(2,N)$ partition \cite{etayo2021spherical}.

The first paper \cite{beltran2020diamond} concentrates on logarithmic energy and cites Zhou's construction. It would be interesting to compare the results for logarithmic energy with the empirical logarithmic energy of the $\EQCode(2,N)$ codes, especially considering that the $\EQCode(2,N)$ codes have a non-random rotation offset on each parallel that maximizes the distance between the code points on adjacent parallels \cite[Section 4.1.2]{leopardi2007distributing}.

The second paper \cite{etayo2021spherical} examines spherical cap discrepancy and also cites the PhD thesis \cite{leopardi2007distributing}, and the Matlab toolbox \cite{leopardi2017recursive}. It is interesting to compare the proof of \cite[Theorem 1.6]{etayo2021spherical} on the upper bound for spherical cap discrepancy of the diamond ensemble with the proof of \cite[Theorem 5.4.1]{leopardi2007distributing} on the same topic for the $\EQCode(d,N)$ codes. As expected, the order of the bound in both proofs coincides for $d=2$. The latter proof involves general $d > 1$ rather than just $d=2$, but it uses a similar argument about the number of regions of an equal area partition that contain the boundary of a spherical cap. See, for example, \cite[Figure 3]{etayo2021spherical}.

\section{Some applications}\label{sec:applications}
\subsection*{Biology and medicine}
\begin{figure}[!ht]
\begin{center}
\includegraphics[angle=0,width=100mm]{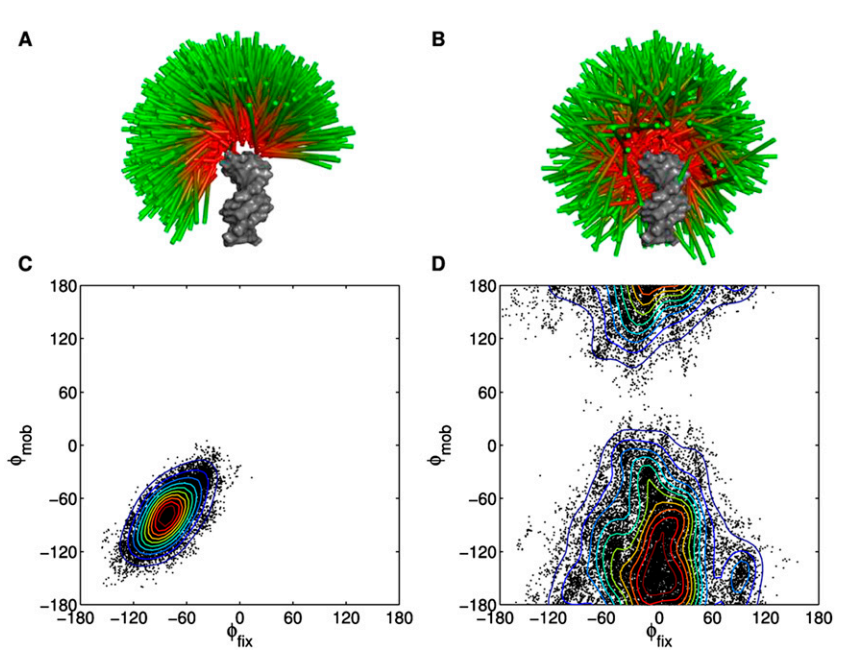}
\caption{
Chu et al.
``Visualization of 1000 randomly selected conformers observed in the dPEG (A) and sPEG (B) simulations,''
\cite{chu2009conformational}.}
\end{center}
\label{fig:chu2009conformational}
\end{figure}

\begin{figure}[!ht]
\begin{center}
\includegraphics[angle=0,width=80mm]{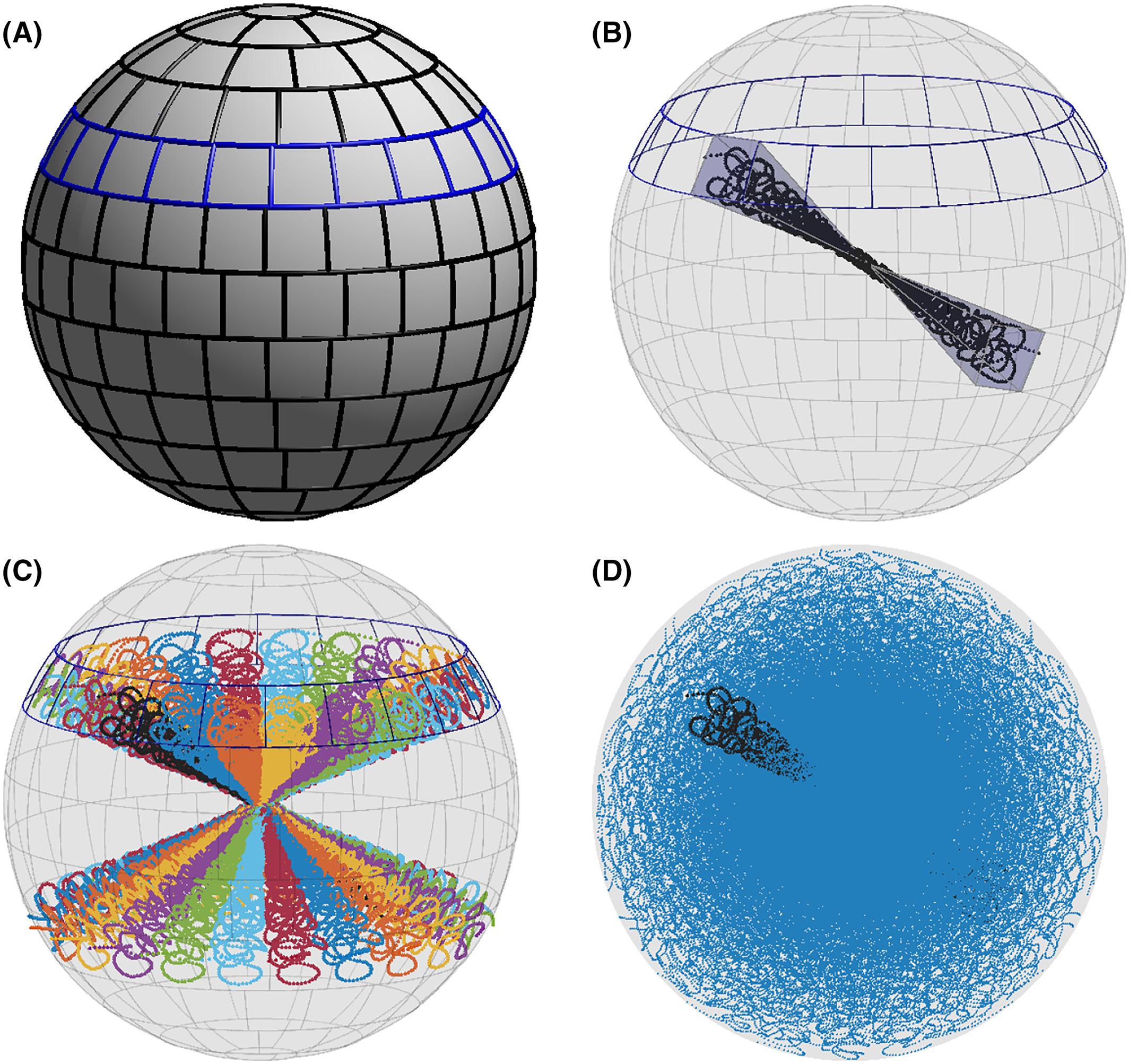}
\caption{
Lazarus et al.
``3D SPARKLING process.
A, Partition of the sphere into 100 regions of equal area.
Regions along a constant elevation angle were highlighted in blue: they are identical up to a rotation.
B, One 3D density sector containing a SPARKLING shot.
C, The SPARKLING shot is then rotated along the considered latitude.
D, the whole fully 3D SPARKLING trajectory. An individual segment is highlighted in black. \dots''
\cite{lazarus20203d}.}
\end{center}
\label{fig:lazarus20203d}
\end{figure}

\begin{figure}[!ht]
\begin{center}
\includegraphics[angle=0,width=100mm]{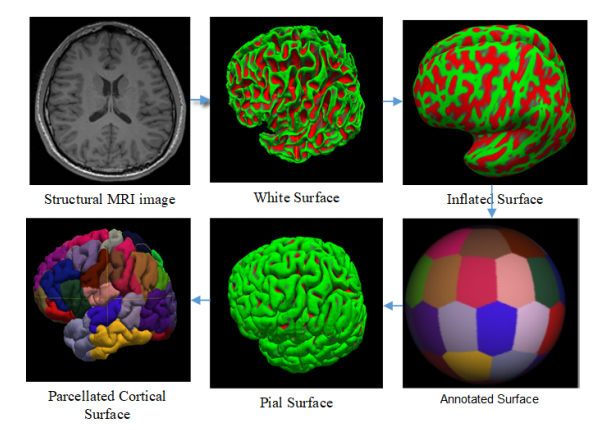}
\caption{
Das and Maharatna,
``Raw MRI scans to parcellated segmented brain image,''
\cite{das2020automated}.}
\end{center}
\label{fig:das2020automated}
\end{figure}

\paragraph*{Biochemistry.}
The 2009 paper by Chu and colleagues \cite{chu2009conformational}
investigates RNA folding by simulating two simple cases where two helices are joined by a non-helix segment.
The methods used in the paper include apparently using $\EQCode(3, 16000)$ to produce
``8000 equally spaced points on the upper half-sphere of the unit three-sphere $\Sphere^3 \subset \Real^4$, yielding a set of quaternions that sampled the space of rigid body rotations $SO(3)$ evenly.'' (Figure \ref{fig:chu2009conformational}).
%
\paragraph*{Medical imaging.}
The 2020 paper by Lazarus and colleagues \cite{lazarus20203d}
extends the ``SPARKLING (Spreading Projection Algorithm for Rapid K-space sampLING)''
optimization algorithm for efficient compressive sampling patterns for 3D magnetic resonance imaging (MRI).
The 3D SPARKLING process uses the $\Partition(2, 100)$ partition to arrange MRI shots in a trajectory.
The paper compares this process with two stacked SPARKLING processes and finds it to be inferior to a variable density stacked
SPARKLING process (Figure \ref{fig:lazarus20203d}).
\paragraph*{Neurobiology.}
The 2020 paper by Das and Maharatna \cite{das2020automated}
presents an ``end-to-end toolchain that processes raw MRI data and generates network metrics for brain connectivity analysis using non-anatomical equal-area parcellation.''
The method presented in this paper is quite involved, but includes steps that use the $\Partition(2,80)$ partition
and $\EQCode(2,80)$ spherical code: ``\ldots we partition the spherical surface into equal sized areas by
applying the equipartition algorithm of unit sphere \cite{leopardi2006partition}. We create a list of centre points of all the equal partitioned areas of unit sphere and scale them up to spherical surface \ldots'' (Figure \ref{fig:das2020automated}).
\subsection*{Climate and weather}
\paragraph*{Climate science.}
The 2018 paper by Werner and colleagues \cite{werner2018spatio}
presents ``the first spatially resolved and millennium-length summer (June–August)
temperature reconstruction over the Arctic and sub-Arctic domain (north of $60^{\circ}$ N).''

\begin{figure}[!ht]
\begin{center}
\includegraphics[angle=0,width=100mm]{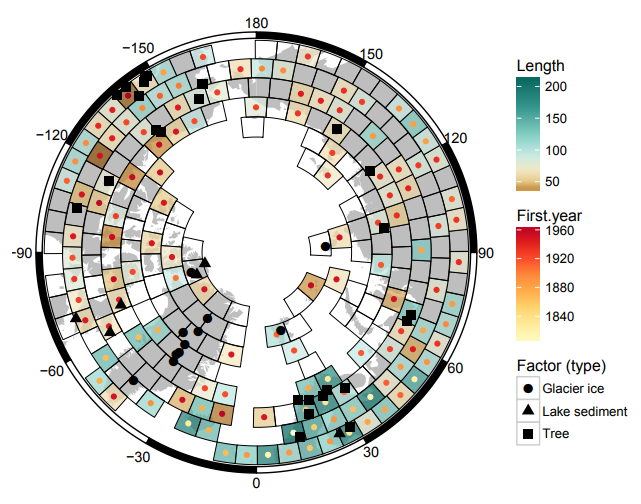}
\caption{
Werner et al.
``Distribution of input data. Length (fill of quadrilaterals)
and first year (coloured circles) of the regridded instrumental data.
Symbols show the locations and type of proxy data used (PAGES 2k
Consortium, 2017). The reconstruction target area is all grid cells
marked with wire frames.''.
\cite{werner2018spatio}.}
\end{center}
\label{fig:werner2018spatio}
\end{figure}

\begin{figure}[!ht]
\begin{center}
\includegraphics[angle=0,width=140mm]{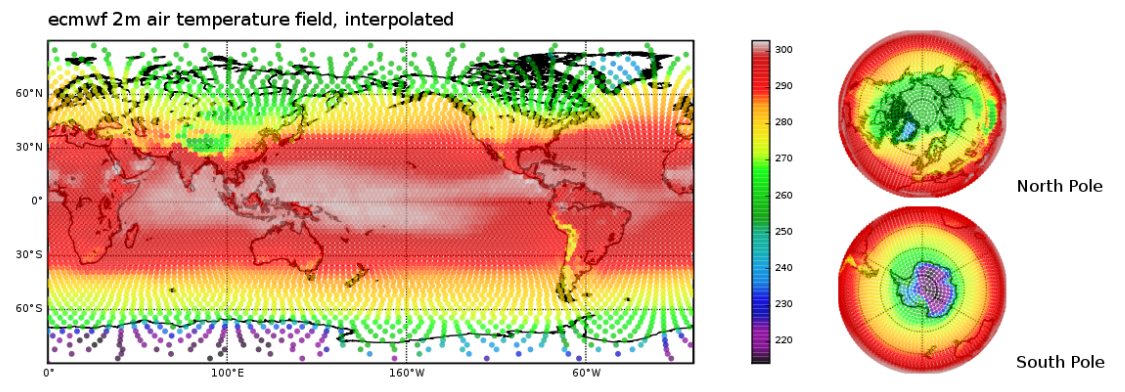}
\caption{
Fauchereau et al.
``ERA 15 surface temperature long-term mean (1979–1993): interpolated onto a zonal equal area partitioning of the sphere using 6500 points.''
\cite{fauchereau2008empirical}.}
\end{center}
\label{fig:fauchereau2008empirical}
\end{figure}

\begin{figure}[!ht]
\begin{center}
\includegraphics[angle=0,width=100mm]{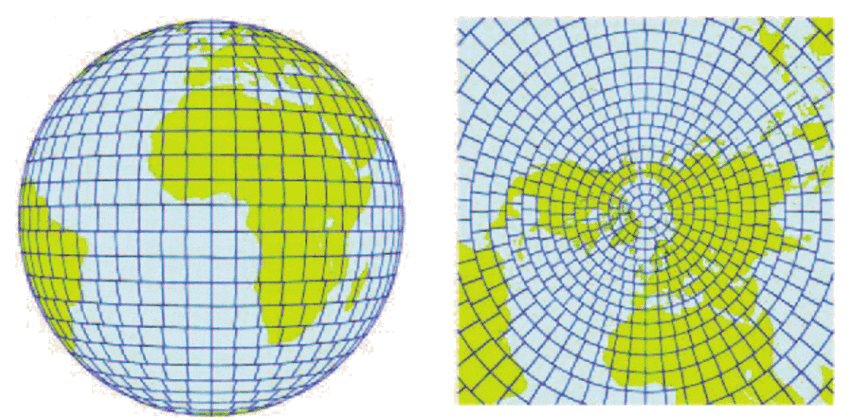}
\caption{Mozdzynski et al.
``EQ\_REGIONS partitioning of grid-point space,''
\cite{mozdzynski2012pgas}.}
\end{center}
\label{fig:mozdzynski2012pgas}
\end{figure}

The 2008 paper by Fauchereau and colleagues \cite{fauchereau2008empirical}
applies Empirical Mode Decomposition (EMD) ``in two dimensions over the sphere to demonstrate its
potential as a data-adaptive method of separating the different scales of spatial variability in a geophysical (climatological/meteorological) field.''
The paper uses the $\Partition(2, 6500)$ partition and the $\EQCode(2, 6500)$ spherical code (Figure \ref{fig:fauchereau2008empirical}).
\paragraph*{Numerical weather prediction.}
Papers by Mozdzynski and others at the European Centre for Medium Range Weather Forecasts (ECMWF) \cite{deconinck2016accelerating,deconinck2017atlas,mozdzynski2012pgas,mozdzynski2015partitioned,wedi2015modelling}
describe the use of code derived from the EQSP Matlab Toolbox \cite{leopardi2005recursive} to balance the parallel
load of the ECMWF Integrated Forecasting System (IFS). The papers call this load balancing method \emph{EQ\_REGIONS partitioning} (Figure \ref{fig:mozdzynski2012pgas}).
\subsection*{Geology and geophysics}
\begin{figure}[!ht]
\begin{center}
\includegraphics[angle=0,width=100mm]{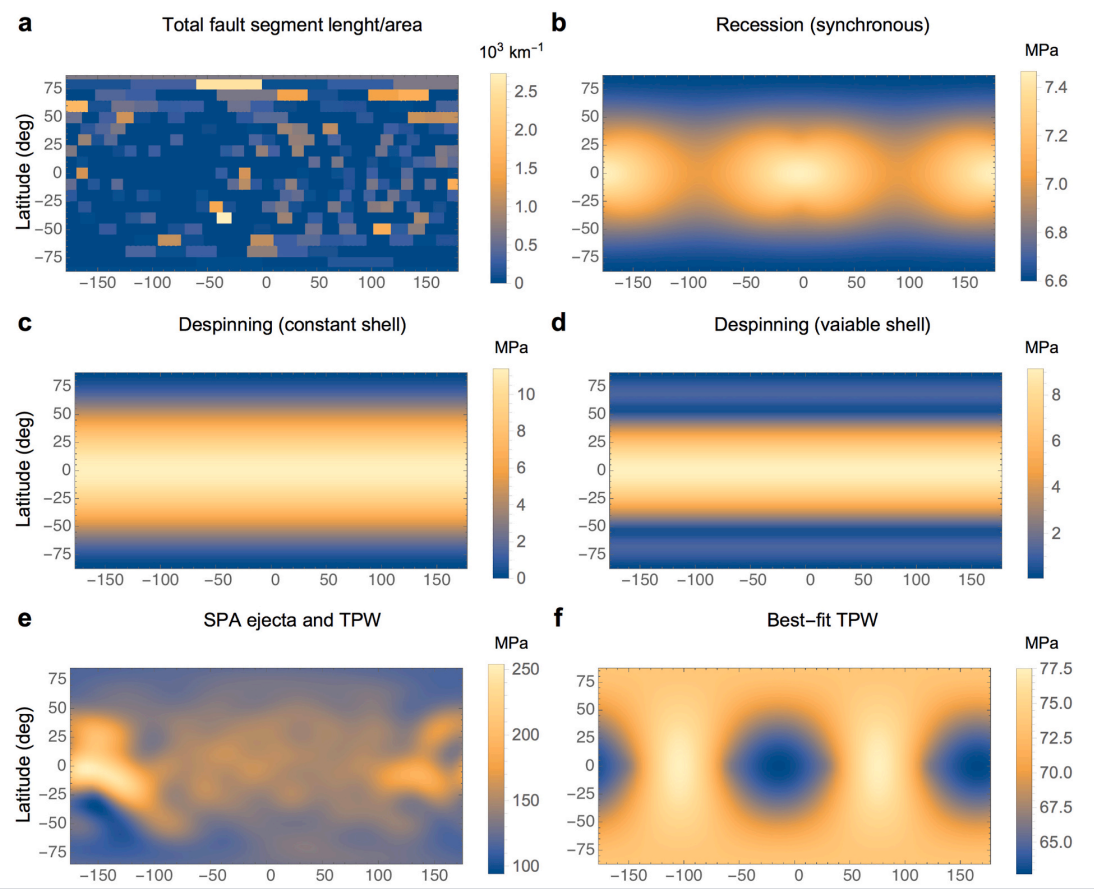}
\caption{
Matsuyama et al.
``a, Total fault segment length per unit area computed by sampling the digitized fault segments in 400 equal area regions partitioned using the `igloo' method of Leopardi (2006).
b-f, Absolute difference between the maximum and minimum principal stresses (principal stresses difference), which quantifies the deviatoric stress, for a variety of mechanisms combined with isotropic contraction \ldots''
\cite{matsuyama2021global}.}
\end{center}
\label{fig:matsuyama2021global}
\end{figure}

Papers and theses describing applications in geology and geophysics include \cite{alken2021evaluation,devriese2014enhanced,devriese2016detecting,domingos2018geomagnetic,hammer2018local,hammer2021applications,hammer2021geomagnetic,hammer2022secular,istas2023transient,kloss2019time,matsuyama2021global,olsen2017lcs,olsen2023determination}.
For example, the paper  by Matsuyama and colleagues \cite{matsuyama2021global}
uses the $\Partition(2, 400)$ partition to sample tectonic patterns on the moon (Figure \ref{fig:matsuyama2021global}).
The paper does not justify this choice of sampling method or estimate its accuracy.

The paper by Alken and colleagues \cite{alken2021evaluation}
uses a robust Huber model based on $10\ 000$ points obtained via the $\EQCode(2, 10\ 000)$ spherical code as a component of the evaluation of models of the Earth's magnetic field.

\begin{figure}[!ht]
\begin{center}
\includegraphics[angle=0,width=100mm]{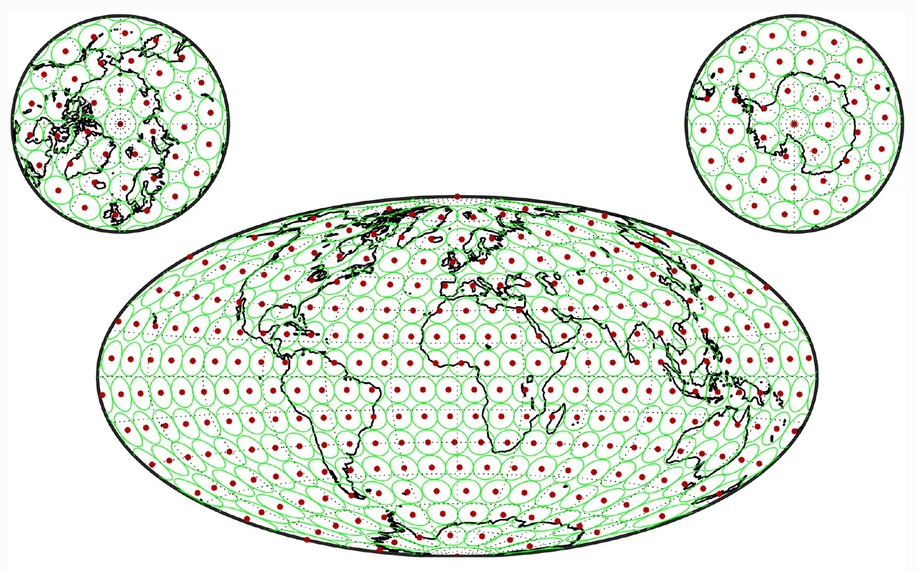}
\caption{
Hammer et al.
``Distribution of the 300 GVOs (red dots) and associated cylinder bins (in green) using a Hammer projection,''
\cite{hammer2021geomagnetic}.}
\end{center}
\label{fig:hammer2021geomagnetic}
\end{figure}

The thesis by Domingos \cite{domingos2018geomagnetic},
and the papers by
Hammer and colleagues \cite{hammer2021applications,hammer2021geomagnetic,hammer2022secular}
Istas and colleagues \cite{istas2023transient},
and Kloss and colleagues \cite{kloss2019time}
use the $\EQCode(2, 300)$ or $\EQCode(2, 500)$ spherical codes to locate either 300 points or 500 points around the Earth,
and use these points to locate geomagnetic virtual observatories (Figure \ref{fig:hammer2021geomagnetic}).
Each of these is effectively an approximate solution of the spherical cap packing problem.

\subsection*{Materials science}
\begin{figure}[!ht]
\begin{center}
\includegraphics[angle=0,width=80mm]{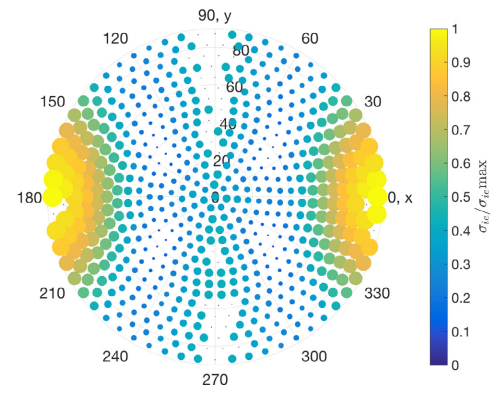}
\caption{
Sabiston et al.
``Interface stress distribution over the fibre ODF for plaque location 1
considering the entire microstructure, excluding orientations that do not
appear in the evaluated microstructural image,''
\cite{sabiston2021accounting}.}
\end{center}
\label{fig:sabiston2021accounting}
\end{figure}

The 2021 paper by Sabiston and colleagues \cite{sabiston2021accounting}
presents and evaluates a micromechanics model for use in the
fatigue characterization of injection moulded carbon fibre.
The microstructure is characterized in terms of the orientation of carbon fibres,
as an orientation distribution function (ODF).
This function is approximated through the use of the $\Partition(2, 1200)$ partition
and the $\EQCode(2, 1200)$ spherical code.
``1200 was selected by performing a parametric study on the effect of number of orientations on
the homogenized stress as well as the maximum interface stress''
(Figure \ref{fig:sabiston2021accounting}).
\subsection*{Mathematical physics}
\begin{figure}[!ht]
\begin{center}
\includegraphics[angle=0,width=120mm]{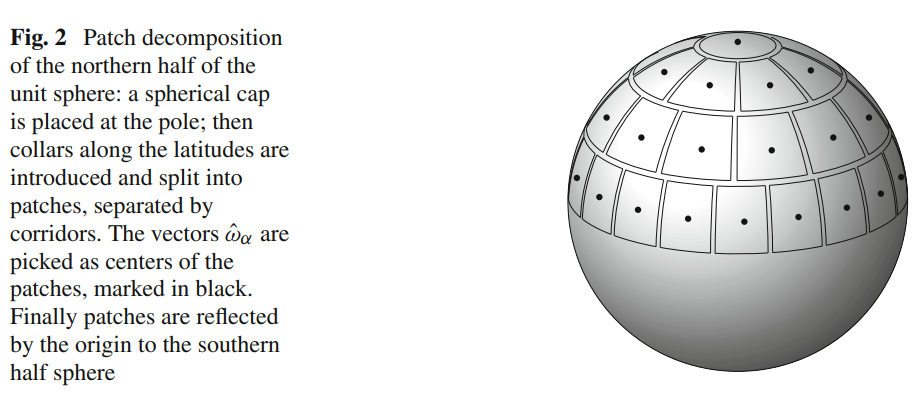}
\caption{
Benedikter et al.
\cite{benedikter2021correlation}.}
\end{center}
\label{fig:benedikter2021correlation}
\end{figure}

The 2021 paper by Benedikter and colleagues \cite{benedikter2021correlation}
rigorously derives ``the leading order of the correlation energy of a
Fermi gas in a scaling regime of high density and weak interaction.''
The paper uses a modified version of the $\Partition(2, M)$ partition that (1)
partitions the northern hemisphere and reflects this partition into the southern hemisphere; and
(2) introduces corridors between the regions;
(Figure \ref{fig:benedikter2021correlation}).

\subsection*{Robotics}
\begin{figure}[!ht]
\begin{center}
\includegraphics[angle=0,width=120mm]{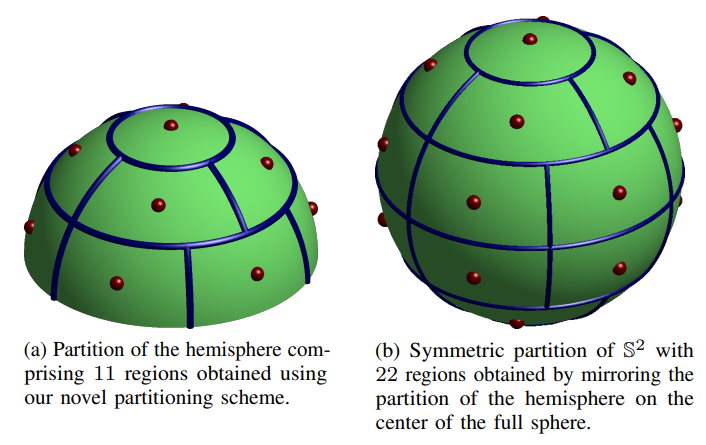}
\caption{
Pfaff et al.
``Illustration showing a partition of $\mathbb{H}^2$ with equally sized regions and
a partition of $\mathbb{S}^2$ obtained by mirroring the partition of the hemisphere,''
\cite{pfaff2020hyperhemispherical}.}
\end{center}
\label{fig:pfaff2020hyperhemispherical}
\end{figure}

The 2020 paper by Pfaff and colleagues \cite{pfaff2020hyperhemispherical}.
proposes ``a grid filter
for arbitrary-dimensional unit hyperhemispheres and apply it to an orientation estimation task and another evaluation scenario.''
(Figure \ref{fig:pfaff2020hyperhemispherical}).
It is one of a series of related papers \cite{frisch2023deterministic,kurz2017discretization,li2020grid,li2020nonlinear,li2021hyperspherical,li2021progressive,pfaff2020hyperhemispherical,pfaff2020spherical,pfaff2022state} that each cite \cite{leopardi2006partition}.
The partitions used in the paper differ from $\Partition(d, N)$ partitions in the following way:
``we adjusted the algorithm so that it
yields the best even integer number of collars. Then, when
subdividing from top to bottom, the boundary of one collar will
run along the equator of the hypersphere.''
The paper does not explain how the modified algorithm generates regions of equal area,
nor does it provide source code.
\subsection*{Visualization}
\begin{figure}[!ht]
\begin{center}
\includegraphics[angle=0,width=100mm]{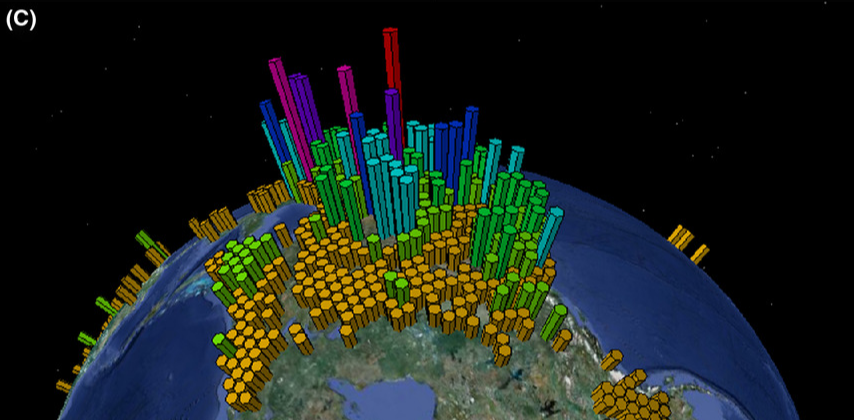}
\caption{Arrigo, et al.
``Species diversity of the \emph{Selaginella} subgenus \emph{Tetragonostachys} in North America using a grid of 20 000 cells with equal areas,''
\cite{arrigo2012quantitative}.}
\end{center}
\label{fig:arrigo2012quantitative}
\end{figure}

The 2012 paper by Arrigo and colleagues \cite{arrigo2012quantitative}
describes the R2G2 R CRAN package for the visualization of spatial data using Google Earth.
The package uses the $\Partition(2, 50)$, $Partition(2, 500)$, $\Partition(2, 5000)$, $\Partition(2, 10\ 000)$, and $\Partition(2, 20\ 000)$
partitions to calculate and plot histograms and other visualizations of data distributed on the Earth's surface
(Figure \ref{fig:arrigo2012quantitative}).
\subsection*{Sphere vs hypersphere}
The majority of application papers citing \cite{leopardi2006partition} focus only on applications on the sphere $\Sphere^2$.
Even so, only about 14 papers citing \cite{leopardi2006partition} also cite either the original work on zonal equal area partitioning of $\Sphere^2$ by Rakhmanov, Saff and Zhou \cite{RakSZ94} or Zhou's PhD thesis \cite{Zho95}.
Of these 14 papers (\cite{ahmadia2010parallel,brauchart2008optimal,brauchart2012next,dickstein2019approximation,hardin2016comparison,holhocs2014octahedral,
ishii2011comparison,jenkins2012local,marantis2009comparison,michaels2017node,sun2008spherical,wang2023numerical,womersley2018efficient,xie2009spherical,zotter2009sampling})
only six (\cite{ahmadia2010parallel,ishii2011comparison,jenkins2012local,marantis2009comparison,xie2009spherical,zotter2009sampling})
are applications-oriented.

Ahmadia's 2010 PhD thesis uses the $\Partition(d,N)$ partitions with $d$ from 5 to 12 to solve a semiconductor lithography optimization problem \cite{ahmadia2010parallel}.
The 2011 paper by Ishii uses the $\Partition(2,N)$ partitions and $\EQCode(2,N)$ codes as possible solutions to a sampling problem used to calculate the total radiated power from radio equipment -- a quadrature problem \cite{ishii2011comparison}).
Jenkins' 2012 paper \cite{jenkins2012local} concerns the construction of sparse spanners of unit ball graphs in $\Real^3$, relating this problem to the covering problem on the sphere, and using $\Partition(2,N)$ partitions to efficiently approximate coverings.
The 2009 paper by Marantis and colleagues \cite{marantis2009comparison} is described in Section \ref{sec:evaluations}.
The 2009 paper by Xie and colleagues \cite{xie2009spherical} is a conference paper companion to \cite{marantis2009comparison}.
The 2009 paper by Zotter \cite{zotter2009analysis} compares a number of different methods of approximating polynomial functions on the sphere $\Sphere^2$ in order to analyze discrete spherical microphone and loudspeaker arrays.
Unfortunately this paper seems to confuse interpolation using extremal fundamental systems on $\Sphere^2$ with hyperinterpolation \cite{sloan2000constructive}.

Of the papers that tackle applications on higher-dimensional spheres, some use the double covering
of the $SO(3)$ group of rotations in $\Real^3$ by the $SU(2)$ group, represented by the unit quaternions,
modelled as the hypersphere $\Sphere^3$, and therefore use the $\Partition(3,N)$ partitions and the $\EQCode(3,N)$ codes \cite{chu2010probing,chu2009conformational,del2021multiscale,fusiello2022exact,mehta2009daisy,mehta2012pegus,pfaff2020hyperhemispherical,pfaff2022state}. A few, such as the 2022 paper by Ram{\'\i}rez and Elvingson \cite{ramirez2022efficient} address $\Sphere^3$ and the $\EQCode(3,N)$ codes for other reasons.
Others address applications in higher dimensions, including:
 Ahmadia's 2010 PhD thesis \cite{ahmadia2010parallel}, as described above;
 a 2012 report by Kessler and colleagues \cite{kessler2012construction} that describes an algorithm that uses the $\EQCode$ codes to construct an approximately optimal path to extinction in systems of arbitrary dimensions;
 the 2016 paper by Rachsinger and colleagues \cite{rachinger2016phase}, as described in Section \ref{sec:evaluations};
 the related 2016 paper by Sedaghat and colleagues \cite{sedaghat2016continuous}, also described in Section \ref{sec:evaluations};
 a 2017 paper by Kurz and Hanebeck \cite{kurz2017linear} that uses the $\EQCode$ codes to construct linear regression Kalman filters; and finally
 a 2021 paper by Miyamoto and colleagues \cite{miyamoto2021constructive} that uses Hopf fibrations to construct
 spherical codes in $\Real^{2^k}$, comparing these to $\EQCode(2^k - 1,N)$ codes for $k$ from 2 to 5.

Of the mathematical papers that use the $\Partition$ partitions and the $\EQCode$ codes in higher dimensions, one stands out:  the 2019 paper by Kunc and Fritzen \cite{kunc2019generation}, which uses the $\EQCode$ codes as starting points
for energy minimization.
\section{Conclusion}
Judging from the wide variety of applications of the $\Partition$ partitions and the $\EQCode$ codes,
these constructions appear to be widely applicable.

Closer inspection reveals that the constructions perform poorly on some problems.
Chief among these is the reconstruction of functions via spherical harmonics \cite{marantis2009comparison},
accomplished on $\Sphere^2$, for example, by scattered data approximation \cite{JetSW98b,JetSW99,sloan2000constructive}.
The joint paper with Sommariva and Vianello \cite{leopardi2017optimal} addresses this problem by using large $\EQCode$ codes as norming sets,
and constructing subsets that have approximation properties almost as good as the norming sets.

In the case covered in this paper,
the applicability of the mathematical construction appears to be not so unreasonable,
given the work that has been done in testing, for each relevant problem,
the performance of the construction relative to alternatives,
and the fitting of the construction into an overall solution that addresses each specific application.

\section*{Acknowledgments}
This paper is based on a series of presentations given at
Oak Ridge National Laboratory (ORNL) in 2014,
The European Centre for Medium-Range Weather Forecasts (ECMWF) in 2016,
The Bureau of Meteorology in 2017, and
The Queensland Association of Mathematics Teachers (QAMT) State Conference in 2021.

Thanks to Kate Evans of ORNL, George Mozdzynski of ECMWF, and Monique Russell of QAMT for invitations to speak.
Thanks to the Bureau of Meteorology for support to visit ECMWF in Reading in 2016.

The original research in the PhD thesis \cite{leopardi2007distributing}
was supported by a University Postgraduate Award from the University of New South Wales.

\section*{Declaration of interests}
The author declares that he has no known competing financial interests or personal relationships that could have appeared to influence the work reported in this paper.

\bibliography{Leopardi-2024-JAS}

\end{document}